\documentclass[12pt]{article}
\usepackage{amsmath,amsfonts,amsthm,amscd,upref,amstext}
\usepackage[dvips]{graphicx}
\usepackage{subfigure}
\usepackage{pb-diagram,pb-xy}
\usepackage[all]{xy}

\newtheorem{prop}{Proposition}[section]
\newtheorem{thm}[prop]{Theorem}

\newtheorem{defn}[prop]{Definition}
\newtheorem{rem}[prop]{Remark}
\newtheorem{lem}[prop]{Lemma}

\newcommand{\N}{\mathbb{N}}
\newcommand{\Z}{\mathbb{Z}}

\begin{document}

\title{A numerical characterization of reduction for arbitrary modules}

\author{R. Callejas-Bedregal$^{1,}\,$\thanks{This work was prepared during
the first author stay at ICMC-USP-S\~ao Carlos-Brazil financed by
CNPq-Brazil - Grant 151733/2006-6.  2000 Mathematics Subject
Classification: 13H15(primary), 13B22(secondary). {\it Key words}:
Rees's theorem, integral closure, multiplicity sequence
}\,\,\,\,and\,\,\,V. H. Jorge P\'erez$^{2}$}

\date{}
\maketitle

\noindent $^1$ Universidade Federal da Para\'\i ba-DM, 58.051-900,
Jo\~ao Pessoa, PB, Brazil ({\it e-mail: roberto@mat.ufpb.br}).

\vspace{0.3cm} \noindent $^2$ Universidade de S{\~a}o Paulo -
ICMC, Caixa Postal 668, 13560-970, S{\~a}o Carlos-SP, Brazil ({\it
e-mail: vhjperez@icmc.usp.br}).

\vspace{0.3cm}
\begin{abstract}

Let $(R, \mathfrak m)$ be a $d$-dimensional Noetherian local ring
and $E$ a finitely generated $R$-submodule of a free module $R^p.$
In this work we introduce a multiplicity sequence $c_k(E),
\,k=0,\ldots, d+p-1$ for $E$ that generalize the Buchsbaum-Rim
multiplicity defined when $E$ has finite colength in $R^p$ as well
as the Achilles-Manaresi multiplicity sequence that applies when
$E\subseteq R$ is an ideal. Our main result is that the new
multiplicity sequence can indeed be used to detect integral
dependence of modules. Our proof is self-contained and implies
known numerical criteria for integral dependence of ideals and
modules.
\end{abstract}

\maketitle

\section{Introduction}

\vspace{0.3cm}

Let $(R, \mathfrak{m})$ be a local Noetherian ring, $N$ a finitely
generated $d$-dimensional $R$-module, and $I\subseteq J$ be two
ideals in $A$. Recall that $I$ is a reduction of $(J,N)$ if
$IJ^nN=J^{n+1}N$ for sufficiently large $n.$  If $I\subseteq J$
are $\mathfrak{m}$-primary and $I$ is a reduction of  $(J,N)$ then
it is well known and easy to prove that the Hilbert-Samuel
multiplicities $e(J,N)$ and $e(I,N)$ are equal. D. Rees proved his
famous result, which nowadays has his name, that the converse also
holds:

\begin{thm}\label{Rees}(Rees's Theorem, \cite{Rees1}) Let $(R, \mathfrak m)$
be a quasi-unmixed local ring, $N$ a finitely generated
$d$-dimensional $R$-module and  $I\subseteq J$  $\mathfrak
m$-primary ideals of $R.$ Then, the following conditions are
equivalent:
\begin{enumerate}
\item [(i)] $I$ is a reduction of $(J,N)$;

\item [(ii)] $e(J,N)=e(I,N).$
\end{enumerate}
\end{thm}

Now assume that $I\subseteq J$ are arbitrary ideals with the same
radicals. If $I$ is a reduction of $J$ then we have always
$e(J_{\mathfrak p},R_{\mathfrak p})=e(I_{\mathfrak p},R_{\mathfrak
p})$ for all minimal primes of $J.$ However, the converse is not
true, in general. Under additional assumption E. B\"{o}ger
\cite{Boger} was able to prove a converse as follows: {\it let
$J\subseteq I\subseteq\sqrt{I}$ be ideals in a quasi-unmixed local
ring $R$ such that $s(I)=\mbox{ht}(I),$ where $s(I)$ denotes the
analytic spread of $I$. Then $I$ is a reduction of $J$ if and only
if $e(J_{\mathfrak p},R_{\mathfrak p})=e(I_{\mathfrak
p},R_{\mathfrak p})$ for all minimal primes of $J.$}

Using the $j$-multiplicity defined by R. Achilles and M. Manaresi
\cite{Achilles-Manaresi1} (a generalization of the classical
Hilbert-Samuel multiplicity), H. Flenner and M. Manaresi
\cite{Flenner-Manaresi} gave numerical characterization of
reduction ideals which generalize B\"{o}ger's theorem to arbitrary
ideals: {\it let $I\subseteq J$ be ideals in a quasi-unmixed local
ring $R$ and $N$ a finitely generated $d$-dimensional $R$-module.
Then $I$ is a reduction of $(J,N)$ if and only if $j(J_{\mathfrak
p},N_{\mathfrak p})=j(I_{\mathfrak p},N_{\mathfrak p})$ for all
${\mathfrak p}\in \mbox{Spec}(R).$}

There is another generalization of the classical Hilbert-Samuel
multiplicity for arbitrary ideals due to R. Achilles and M.
Manaresi \cite{Achilles-Manaresi}. They introduced, for each ideal
$I$ of a $d$-dimensional local ring $(R, \mathfrak{m})$ and $N$ a
finitely generated $d$-dimensional $R$-module, a sequence of
multiplicities $c_0(I,N),\dots ,c_d(I,N)$ which generalize the
Hilbert-Samuel multiplicity in the sense that for
$\mathfrak{m}-$primary ideals $I$, $c_0(I,N)$ is the
Hilbert-Samuel multiplicity of $I$ in $N$ and the remaining
$c_k(I,N),\;k=1,\ldots,d$ are zero. In fact, their definition was
given in the case that $N=R$ but their construction can be readily
extended to this context.

Using the above multiplicity sequence defined by R. Achilles and
M. Manaresi  \cite{Achilles-Manaresi}, the authors gave the
following numerical characterization of reduction of ideals which
generalize Rees's theorem for arbitrary ideals \cite[Theorem
5.5]{Bedregal-Perez2}:

\begin{thm}
Let $(R, \mathfrak{m})$ be a quasi-unmixed $d$-dimensional local
ring and $N$ a finitely generated $d$-dimensional $R$-module. Let
$I\subseteq J$ be proper arbitrary ideals of $R$ such that
$\mbox{ht}_N(I)>0.$ Then the following conditions are equivalent:

\begin{enumerate}
\item [(i)] $I$ is a reduction of $(J,N)$;

\item [(ii)] $c_k(I,N)=c_k(J,N)$ for all $k=0,...,d$.
\end{enumerate}
\end{thm}

On the other hand, the Buchsbaum-Rim multiplicity $e_{BR}(E)$ is a
generalization of the Samuel multiplicity and is defined for
submodules of free modules $E\subset R^p$ such that $R^p/E$ has
finite length. These were first described by D. A. Buchsbaum and
D. S. Rim in \cite{Buchsbaum-Rim}. The Buchsbaum-Rim multiplicity
has been generalized, in the finite colength case, by D. Kirby
\cite{Kirby}, D. Kirby and D. Rees \cite{Kirby-Rees}, D. Katz
\cite{Katz}, S. Kleiman and A. Thorup \cite{Kleiman-Thorup} and A.
Simis, B. Ulrich and W.Vasconcelos
\cite{Simis-Ulrich-Vasconcelos}. For an extensive history of
Buchsbaum-Rim multiplicity we refer to \cite{Kleiman-Thorup}.
Using the Buchsbaum-Rim multiplicity S.Kleiman and A. Thorup
\cite{Kleiman-Thorup}, D. Katz \cite{Katz} and A. Simis, B. Ulrich
and W.Vasconcelos \cite{Simis-Ulrich-Vasconcelos} proved the
following generalization of the Rees's theorem for modules:

\begin{thm}\label{Rees-mod-finite} Let $(R, \mathfrak m)$
be a quasi-unmixed local ring, $E\subseteq F$  finitely generated
$R$-submodule of a free module $R^p$ such that $R^p/E$ has finite
length. Then, the following conditions are equivalent:
\begin{enumerate}
\item [(i)] $E$ is a reduction of $F$;

\item [(ii)] $e_{BR}(E)=e_{BR}(F).$
\end{enumerate}
\end{thm}

In the last fifteen years the Buchsbaum-Rim  multiplicity of a
submodule of a free module has played an important role in the
theory of equisingularity of families of complete intersections
with isolated singularities (ICIS). The Buchsbaum-Rim multiplicity
has been used in the context to control the $A_f,\;W_f, A$ and $W$
conditions of equisingularity, which are analogous to the Whitney
conditions (cf. \cite{Gaffney1}, \cite{Gaffney-Kleiman},
\cite{Gaffney-Massey}, \cite{Kleiman-Thorup2} and the reference
therein). The usefulness of the Buchsbaum-Rim multiplicity is
restricted to families of ICIS, because it is only for these
singularities that the submodules associated to the
equisingularity conditions have finite colength and only for these
types is the Buchsbaum-Rim multiplicity well defined. In order to
generalize those works for families of arbitrary complete
intersection singularities (ACIS) it is strictly necessary to
generalize first the notion of Buchsbaum-Rim multiplicities for
submodules $M$ of a free module $F$ of arbitrary colength. For
this new notion of multiplicity to be useful in equisingularity
theory it must characterize the integral closure of arbitrary
modules, that is, it must generalize Theorem
\ref{Rees-mod-finite}.

There have been some generalizations of the Buchsbaum-Rim
multiplicity for arbitrary submodules $E$ of a free module $R^p$
which we now describe. T. Gaffney in \cite{Gaffney} introduced a
sequence of multiplicities $e_i(E),\;0\leq i\leq d=\dim R$ in the
analytic context. This sequence satisfies a Rees type theorem:
Suppose that $E\subset F\subset R^p$ are $R:={\cal
O}_{X,x}$-modules where $X^d$ is a complex analytic space which as
a reduced space is equidimensional, and which is generically
reduced. Suppose that $e_i(E, x)=e_i(F, x),\;0\leq i\leq d.$ Then
$E$ is a reduction os $F$. Also, if $E$ is of finite colength in
$R^p$, then $e_d(E)$ is the standard Buchsbaum-Rim multiplicity of
$E$, and the others $e_i$'s are zero. Unfortunately, for ideals of
non-finite colength, Gaffney's multiplicity sequence does not
coincide with the Achilles-Manaresi multiplicity sequence and also
the codimension condition of $E$ in $R^p$ is built into the
definition of the multiplicity which uses a codimension filtration
ascending from the integral closure of the module.

On the other hand, the authors in \cite{Bedregal-Perez1} extended
the notion of the Buchsbaum-Rim multiplicity of a submodule of a
free module to the case where the submodule no longer has finite
colength. For a submodule $E$ of $R^p$ they introduced a sequence
$e_{BR}^k(E),\;k=0,\cdots,d+p-1$ which in the ideal case coincides
with the multiplicity sequence $c_0(I,R),\ldots,c_d(I,R)$ defined
for an arbitrary ideal $I$ of $R$ by R. Achilles and M. Manaresi
\cite{Achilles-Manaresi}. They also proved that if $E=I_1\oplus
\ldots \oplus I_p\subset R^p$ has finite colength then
$e_{BR}^0(E)=p\;!(e_{BR}(E))$ and $e_{BR}^k(E)=0$ for
$k=1,\ldots,d-1$. Nevertheless, no relation with reduction of
modules and their multiplicity sequence was shown in their work.

There is  also a particularly beautiful generalization of
Flenner-Manaresi theorem for arbitrary submodules of a free
modules due to B. Ulrich and J. Validashti (see
\cite{Ulrich-Validashti}), they introduced a multiplicity $j(E)$
for a submodule of the free module $R^p$ that generalizes the
Buchsbaum-Rim multiplicity defined when $E$ has finite colength in
$R^p$ as well as the $j$-multiplicity of Achilles-Manaresi that
applies when $E\subseteq R$ is an ideal. Their result is as
follows:

\begin{thm}\label{Boger-mod} Let $(R, \mathfrak m)$ be a universally
catenary ring, $E\subseteq F$ finitely generated $R$-submodule of
a free module $R^p,$ and $N$ a finitely generated locally
equidimensional Noetherian $R$-module. Assume that $E_{\mathfrak
p}=F_{\mathfrak p}$ for every minimal prime $\mathfrak p$ of $R.$
Then, the following are equivalent:
\begin{enumerate}
\item [(i)] $E$ is a reduction of $F$;

\item [(ii)] $j(E_{\mathfrak q})=j(F_{\mathfrak q})$ for every
${\mathfrak q}\in \mbox{Spec}(R).$
\end{enumerate}
\end{thm}

The above theorem characterize reduction of arbitrary modules by
using numerical data in all localizations of the modules which is
hard to verify algebraically and doesn't seems (at least for the
authors) to be useful in equisingularity theory. In this work we
introduce a multiplicity sequence $c_k(E,N)$ with $k=0,\ldots,
d+p-1$ for the pair $(E, N)$ that generalize the Buchsbaum-Rim
multiplicity defined when $E$ has finite colength in $R^p$ as well
as the Achilles-Manaresi multiplicity sequence that applies when
$E\subseteq R$ is an ideal. Our main result is that the new
multiplicity sequence can indeed be used to detect integral
dependence of modules:

\begin{thm}
Let $(R, \mathfrak{m})$ be a Noetherian local ring, $E\subseteq
F\subseteq R^p$ be $R$-modules and write $I:={\cal R}_1(E)$ for
the corresponding ideal of $A:=\mbox{Sym}(R^p).$ Let $N$ be a
$d$-dimensional finitely generated $R$-module and set
$M:=A\otimes_R N.$ Assume that $\mbox{ht }_M(I)>0.$ Consider the
following statements:

\begin{enumerate}
\item [(i)] $E$ is a reduction of $(F,N)$;

\item [(ii)] $c_k(E,N)=c_k(F,N)$ for all $k=0,...,d+p-1$.
\end{enumerate}
Then, $(i)$ implies $(ii)$ and if $N$ is quasi-unmixed the
converse also holds.
\end{thm}

We strongly believe that this multiplicity sequence, apart of
being important in commutative algebra, it will be very useful for
studying equisingularity conditions for families of ACIS. In
particular we expect that it will characterize the $A_f$
condition, answering positively Gaffney and Kleiman's conjecture
stated in \cite[p. 546]{Gaffney-Kleiman}. For this generalization
of all equisingularity conditions to be carried out it is also
necessary to develop a geometric theory for this multiplicity
sequence, a theory involving blowups and intersection numbers as
in the work of J. P. Henry and M. Merle \cite{Henry-Merle} and S.
Kleiman and A. Thorup \cite{Kleiman-Thorup}. The authors hope to
present this geometric approach elsewhere.

The paper is organized as follows. In section $2$, we recall the
basic results of Hilbert functions of bigraded algebras and we
define the $c^D$-multiplicity sequence associated to a graded
module. The important result of this section is the additivity
formula for this multiplicity sequence. In sections $3$ and $4,$
we define two multiplicity sequences associated to ideals
generated by linear forms, which we call $c^*$-multiplicity
sequence and $c^{\sharp}$-multiplicity sequence. The important
results of this sections are the additivity formula for those
multiplicity sequences which immediately implies that they remains
constant when passing to a reduction (Theorem \ref{Rees-fraco} and
Theorem \ref{Rees-fraco2}). They are related by a third
multiplicity sequence, which is denoted by $b^D$, which also
satisfies the additivity formula (see Lemma \ref{lemma2}). The
$c^*, b^D$ and $c^{\sharp}$-multiplicity sequences serve different
purposes: the first two are more readily seen to be additive on
short exact sequences of graded modules (Theorem \ref{additivity}
and the proof of Proposition \ref{additivity2}) and they were
introduced in this work with the only purpose of proving the
additivity property of the $c^{\sharp}$-multiplicity sequences
(Proposition \ref{additivity2}). The last multiplicity sequence on
the other hand is more suited for proving that conversely, the
constancy of the multiplicity sequence implies integral dependence
(Theorem \ref{reesthm}). In section $5$, we recall the notion of
intertwining algebras and modules and state the reduction
criterium we use here which was proved for algebras by A. Simis,
B. Ulrich and W. Vasconcelos \cite{Simis-Ulrich-Vasconcelos}. The
important result of this section is Theorem \ref{reesthm} which
contain as special case the multiplicity sequence of an arbitrary
module, defined in section $6$, which in turn generalize the
Buchsbaum-Rim multiplicity defined only for finite colength
modules as well as the Achilles-Manaresi multiplicity sequence
defined for arbitrary ideals. The main result of section $6$ is
Theorem \ref{reesmodule} which is an immediate consequence of
Theorem \ref{reesthm}. Our approach is partly inspired by
\cite{Bedregal-Perez2}, \cite{Ulrich-Validashti} and
\cite{Achilles-Manaresi}.

\section{Multiplicity sequence}

In this section we recall some well-known facts on Hilbert
functions and Hilbert polynomials of bigraded modules, which will
be essential for defining the Buchsbaum-Rim multiplicity sequence
associated to a pair $(I,M).$

Let $R=\oplus_{i,j=0}^{\infty}R_{i,j}$  be a bigraded ring and let
$T=\oplus_{i,j=0}^{\infty}T_{i,j}$ be a bigraded $R$-module.
Assume that $R_{0,0}$ is an Artinian ring and that $R$ is finitely
generated as an $R_{0,0}$-algebra by elements of $R_{1,0}$ and
$R_{0,1}$ (i.e., $R$ is a standard bigraded algebra) The {\it
Hilbert function} of $T$ is defined to be
$$h_T(i,j)=\ell_{R_{0,0}} (T_{i,j}).$$
For $i,j$ sufficiently large, the function $h_T(i,j)$ becomes a
polynomial $P_T(i,j).$ If $D$ denotes the dimension of the module
$T$, we can write this polynomial in the form
$$P_T(i,j)= \sum_{\begin{array}{l}
k,l\geq 0\\
k+l \leq D-2\end{array}}\!\!\!\!\!\!\!\!\!\!\!a_{k,l}(T)\left(
\begin{array}{c}
i+k\\
k\end{array}\right) \left( \begin{array}{c}
j+l\\
l\end{array}\right)$$ \noindent with $a_{k,l}(T)\in \mathbb Z$ and
$a_{k,l}(T)\geq 0$ if $k+l= D-2$ \cite[Theorem 7, p. 757 and
Theorem 11, p. 759]{Waerden}.

We also consider the sum transform of $h_T$ with respect to the
first variable defined by
$$h_T^{(1,0)}(i,j)=
\sum_{u=0}^{i}h_T(u,j).$$ From this description it is clear that,
for $i,j$ sufficiently large, $h_T^{(1,0)}$  becomes a polynomial
with rational coefficients of degree at most $D-1.$ As usual, we
can write this polynomial in terms of binomial coefficients

$$P_T^{(1,0)}(i,j)= \sum_{\begin{array}{l}
k,l\geq 0\\
k+l \leq
D-1\end{array}}\!\!\!\!\!\!\!\!\!\!\!a_{k,l}^{(1,0)}(T)\left(
\begin{array}{c}
i+k\\
k\end{array}\right)\left( \begin{array}{c}
j+l\\
l\end{array}\right)$$ \noindent with $a_{k,l}^{(1,0)}(T)$ integers
and $a_{k,D-k-1}^{(1,0)}(T) \geq 0,.$

Since
$$h_T(i,j)=h_T^{(1,0)}(i,j)-h_T^{(1,0)}(i-1,j)$$
\noindent we get $a_{k+1,l}^{(1,0)}(T)=a_{k,l}(T)$ for $k,l\geq
0,\,k+l\leq D-2.$

\begin{defn}
For the coefficients of the terms of highest degree in
$P_T^{(1,0)}$ we introduce the symbols
$$c_{k}(T):=a_{k,D-k-1}^{(1,0)}(T),\;\;\;\;k=0,\cdots, D-1$$
\noindent which are called the {\bf multiplicity sequence}  of
$T$.
\end{defn}

We define next the $c^D$-multiplicity sequence associated to a
module. Let $(R,\mathfrak{m})$ be a local ring, $S=\oplus_{j\in
\N}S_j$ a standard graded $R$-algebra, $N=\oplus_{j\in \N}N_j$ a
finitely generated graded $S$-module, and
$$T:=G_{\mathfrak m}(N)=\bigoplus_{i,j\in\N}\,\frac{{\mathfrak m}^iN_j}
{{\mathfrak m}^{i+1}N_j}$$ \noindent the bigraded $F$-module with
$$F:=G_{\mathfrak m}(S)=\bigoplus_{i,j\in\N}\,\frac{{\mathfrak m}^iS_j}
{{\mathfrak m}^{i+1}S_j}.$$ \noindent Notice that  $F_{0,
0}=R/\mathfrak{m}$ is a field.

\begin{defn} Consider an integer $D$ such that $D\geq \dim N.$ For all \linebreak
$k=0, \ldots, D-1,$ we set

$$ c_k^D(N)=\left\{\begin{array}{ll}
               0 & \mbox{ if } \dim N<D \\
               c_k(T) & \mbox{ if } \dim N=D

                \end{array}
                \right. $$
\noindent which is called the {\bf $c^D$-multiplicity sequence} of
$N.$ Moreover, we set \linebreak $c_k(N):=c_k^{\dim N}(N).$
\end{defn}

First we show that this $c^D$-multiplicity sequence behaves well
with respect to short exact sequences.

\begin{prop}(\cite[Proposition 2.3]{Bedregal-Perez2})\label{proposition1}
Let $(R, \mathfrak{m})$ be a local ring, $S=\oplus_{j\in \N}S_j$ a
standard graded $R$-algebra, and $0 \longrightarrow
N_0\longrightarrow N_1\longrightarrow N_2\longrightarrow 0$ an
exact sequence of finitely generated graded $S$-modules. Then for
$D\geq d:=\dim N_1$
$$c_k^D(N_1)=c_k^D(N_0)+c_k^D(N_2)$$ for all $k=0,\ldots, D-1.$
\end{prop}

\begin{proof}
Let $M_s:={\mathcal R}(\mathfrak m, N_s)^+:=\oplus_{i\in
\Z}\oplus_{j\in\N}\,{\mathfrak m}^i(N_s)_j$ be the extended Rees
module associated to $N_s,\, s=0,1,2.$ For any bigraded module $T$
and for $i,j\gg 0,$ we define the polynomial $h^D_{T}(i,j)$ of
degree $D-2$ as the Hilbert polynomial of $h_{T}(i,j)$ adding
coefficient zero to the terms of degree between $\dim(T)-2$ and
$D-2.$

Let $u$ be an indeterminate, which we consider with degree one.
Set \linebreak $M'_0:=\ker(M_1\rightarrow
M_2)=\bigoplus_{i\in\Z,\,j\in\N} (N_0)_j\cap {\mathfrak
m}^i(N_1)_j.$ We consider the natural diagram

$$
\xymatrix { 0 \ar[r] & M_0'(1,0) \ar[r] \ar[d]^{u^{-1}} & M_1(1,0)
\ar[r] \ar[d]^{u^{-1}} & M_2(1,0) \ar[r] \ar[d]^{u^{-1}} & 0 \\
0 \ar[r] & M_0' \ar[r] & M_1 \ar[r]  & M_2 \ar[r] & 0 }
$$
\noindent which gives an exact sequence of cokernels

\begin{equation}\label{equation51}
0 \longrightarrow G':=\frac{M'_0}{u^{-1}M'_0}\longrightarrow
G_{\mathfrak m}(N_1)\longrightarrow G_{\mathfrak
m}(N_2)\longrightarrow 0.
\end{equation}

Denote the cokernel of the natural injection $M_o\hookrightarrow
M_0'$ by $L.$ Using the diagram

$$
\xymatrix { 0 \ar[r] & M_0(1,0) \ar[r] \ar[d]^{u^{-1}} & M'_0(1,0)
\ar[r] \ar[d]^{u^{-1}} & L(1,0) \ar[r] \ar[d]^{u^{-1}} & 0 \\
0 \ar[r] & M_0 \ar[r] & M_0' \ar[r]  & L \ar[r] & 0 }
$$
\noindent the snake-lemma yields an exact sequence

\begin{equation}\label{equation52}
0 \longrightarrow V\longrightarrow G_{\mathfrak
m}(N_0)\longrightarrow G' \longrightarrow  W \longrightarrow 0.
\end{equation}
\noindent where $V$ and $W$ are the kernel and cokernel of
$u^{-1}:L(1,0)\to L$ respectively, i.e., we have the exact
sequence

\begin{equation}\label{equation53}
0 \longrightarrow V\longrightarrow L(1,0)\longrightarrow
L\longrightarrow  W \longrightarrow 0.
\end{equation}

For $n\leq 1$ the coefficient modules of $u^n$ in ${\mathcal
R}(\mathfrak m, N_0)^+$ and in $M_0'$ coincide, hence the action
of $u^{-1}$ on $L$ is nilpotent. Therefore the dimension of $L$ is
at most that of $G',$ which is bounded by $D.$ Thus all modules
occurring in the exact sequence (\ref{equation53}) have dimension
at most $D.$

Now (\ref{equation51}), (\ref{equation52}) and (\ref{equation53})
are exact sequences of finitely generated modules of dimension at
most $D.$ We denote by $h_{N_s}(i,j)$ the Hilbert-Samuel function
of $G_{\mathfrak m}(N_s).$

From (\ref{equation51}) and (\ref{equation52})  we have

\begin{equation}\label{equation54}
h^{D\,(1,0)}_{N_0}(i,j)+h^{D\,(1,0)}_{N_2}(i,j)-h^{D\,(1,0)}_{N_1}(i,j)=h^
{D\,(1,0)}_V(i,j)-h^{D\,(1,0)}_W(i,j).
\end{equation}

Because of (\ref{equation53}) we have

\begin{equation}\label{equation55}
h^{D\,(1,0)}_V(i,j)-h^{D\,(1,0)}_W(i,j)=h^{D\,(1,0)}_L(i+1,j)-h^{D\,(1,0)}_L
(i,j)=h^{D}_L(i,j)
\end{equation}
Hence by (\ref{equation54}) and (\ref{equation55})

$$h^{D\,(1,0)}_{N_0}(i,j)+h^{D\,(1,0)}_{N_2}(i,j)-h^{D\,(1,0)}_{N_1}(i,j)=h^
{D}_L(i,j)$$ \noindent is a polynomial of degree at most $D-2,$
which concludes the proof.
\end{proof}

\section{$c^*$-multiplicity sequence}

In this section we introduce the $c^*$-multiplicity sequence. The
main idea here is to consider a suitable grading on the extended
Rees module as in the work of B. Ulrich and J. Validashti
\cite{Ulrich-Validashti}.

Let $(R,\mathfrak{m})$ be a Noetherian local ring, $A$ a standard
graded Noetherian $R$-algebra, $I$ an ideal of $A$ generated by
elements of degree one and $M$ a finitely generated graded
$A$-module.

Let $t$ be a variable. Consider the extended Rees ring of $I$
$${\mathcal R}(I,A)^{+}:=\oplus_{i\in \Z}I^it^i\subseteq A[t,t^{-1}],$$
\noindent and the extended Rees module
$${\mathcal R}(I,M)^{+}:=\oplus_{i\in \Z}I^iMt^i\subseteq M\otimes_R R[t,t^{-
1}],$$ \noindent where we set $I^i=R$ for $i\leq 0.$ Notice that
${\mathcal R}(I,M)^{+}$ is a module over ${\mathcal R}(I,A)^{+}$
which gives rise to the {\it associated graded module} of $M$ with
respect to $I$,

$$G_I(M):=\frac{{\mathcal R}(I,M)^{+}}{t^{-1}{\mathcal R}(I,M)^{+}}=\oplus_{i\in \N}
\frac{I^iM}{I^{i+1}M}\,t^i,$$ \noindent which is a module over the
associated graded ring $G_I(A)$ of the same dimension as $M.$

Assigning degree zero to the variable $t,$ the Laurent polynomial
ring $A[t, t^{-1}]$ becomes a standard graded Noetherian $R[t,
t^{-1}]$-algebra, and \linebreak $M[t, t^{-1}]:=M\otimes_R R[t,
t^{-1}]$ a finitely generated graded module over this algebra. The
extended Rees ring ${\mathcal R}(I,A)^{+}$ is a homogeneous
$R[t^{-1}]$-subalgebra of $A[t, t^{-1}],$ and hence a standard
graded Noetherian $R[t^{-1}]$-algebra. Furthermore ${\mathcal
R}(I,M)^{+}$ is a homogeneous ${\mathcal R}(I,A)^{+}$-submodule of
$M[t, t^{-1}],$ thus a finitely generated graded module over
${\mathcal R}(I,A)^{+}.$ With respect to this grading,
$G_I(A):={\mathcal R}(I,A)^{+}/t^{-1}{\mathcal R}(I,A)^{+}$
becomes a standard graded Noetherian $R$-algebra and
$G_I(M):={\mathcal R}(I,M)^{+}/t^{-1}{\mathcal R}(I,M)^{+}$ a
finitely generated graded module over this algebra. Notice that

$$[G_I(M)]_n=\oplus_{i\in \N}[I^iM/I^{i+1}M]_n.$$
The grading so defined on the extended Rees module and the
associated graded module is called {\it internal grading}-for it
is induced by the grading on the module $M$ (see
\cite{Ulrich-Validashti}).

\begin{defn}
Let $D$ be any integer with $D\geq \dim M.$ We define the
\linebreak {\bf $c^*$-multiplicity sequence} of $M$ with respect
to $I,$ as
$$c^*_{k,\;D}(I,M):=c^D_k(G_I(M)), \;k=0,\dots, D-1.$$
\noindent where $G_I(M)$ is graded by the internal grading. In the
case where $D=\dim M$ we simply write $c^*_k(I,M)$ instead of
$c^*_{k,\;\dim M}(I,M), \;k=0,\dots, \dim M-1.$
\end{defn}

To be more explicit, consider the standard bigraded $R$-algebra
\linebreak
$S^*:=G_{\mathfrak{m}}(G_I(A))=\oplus_{s,n=0}^{\infty}S^*_{s,n}$
with
$$S^*_{s,n}=\oplus_{i=0}^{\infty}\left [\frac{\mathfrak{m}^sI^iA+I^{i+1}A}
{\mathfrak{m}^{s+1}I^iA+I^{i+1}A}\right ]_n,$$ \noindent where
$G_I(A)$ is graded by the internal grading, and the finitely
generated bigraded module over this algebra
$$T^*=G_{\mathfrak{m}}(G_I(M))=\oplus_{s,n=0}^{\infty}T^*_{s,n}$$
\noindent with
$$T^*_{s,n}=\oplus_{i=0}^{\infty}\left [\frac{\mathfrak{m}^sI^iM+I^{i+1}M}
{\mathfrak{m}^{s+1}I^iM+I^{i+1}M}\right ]_n,$$ \noindent where
$G_I(M)$ is graded by the internal grading.

Observe that  $S^*_{0, 0}=R/\mathfrak{m}$ is a field and $T^*$ has
dimension $\dim M$. We denote the Hilbert-Samuel function
$\ell_{S^*_{0,0}}(T^*_{s,n})$ of $T^*=G_{\mathfrak{m}}(G_I(M))$ by
$h^*_{(I,M)}(s,n)$ and its first Hilbert sum by
$h^{*\,(1,0)}_{(I,M)}(s,n).$ Thus

$$h^*_{(I,M)}(s,n)=\sum_{i=0}^{\infty}\ell_R \left [\frac{\mathfrak{m}
^sI^iM+I^{i+1}M}{\mathfrak{m}^{s+1}I^iM+I^{i+1}M}\right ]_n$$
\noindent and
$$h^{*\,(1,0)}_{(I,M)}(s,n)=\sum_{i=0}^{\infty}\ell_R \left [\frac{I^iM}
{\mathfrak{m}^{s+1}I^iM+I^{i+1}M}\right ]_n.$$

For $s,n\gg 0,$ we define the polynomial $h^{*\, D}_{(I,M)}(s,n)$
of degree $D-2$ as the Hilbert polynomial of $h^*_{(I,M)}(s,n)$
adding coefficient zero to the terms of degree between $\dim M-2$
and $D-2.$

Thus, if $s,n\gg 0,$ the sequence $c^*_{k,\; D}(I,M),
\,k=0,\ldots, D-1$ are the numerators of the leading coefficients
of the polynomial $h^{*\, D\, (1,0)}_{(I,M)}(s,n).$

We will need the fact that the $c^*$-multiplicity sequence is
additive on short exact sequences:

\begin{thm}({\it Additivity})\label{additivity}
Let $(R, \mathfrak{m})$ be a local ring, $A$ a standard graded
Noetherian $R$-algebra, and $I$ an ideal of $A$ generated by
linear forms. If \linebreak $0 \longrightarrow M_0\longrightarrow
M_1\longrightarrow M_2\longrightarrow 0$ is an exact sequence of
finitely generated graded $A$-modules and $D$ an integer with
$D\geq d:=\dim M_1.$ Then, for \linebreak $s,n\gg 0,$
$$h_{(I,M_0)}^{*\,D\,(1,0)}(s,n)+h_{(I,M_2)}^{*\,D\,(1,0)}(s,n)-h_{(I,M_1)}^
{*\,D\,(1,0)}(s,n)$$ \noindent is a polynomial of degree at most
$D-2.$ In particular we have that
$$c^*_{k,\;D}(I,M_1)=c^*_{k,\;D}(I,M_0)+c^*_{k,\;D}(I,M_2)$$ for all
$k=0,\ldots, D-1.$
\end{thm}

\begin{proof}
Let $N_j:={\mathcal R}(I,M_j)^+$ be the extended Rees module
associated to $M_j,$ graded by the internal grading. Set
$$N'_0:=\ker(N_1\rightarrow N_2)=\oplus_{n\in\Z}
\oplus_{i=0}^{\infty} [M_0\cap I^iM_1]_n.$$ We consider the
natural diagram

$$
\xymatrix { 0 \ar[r] & N_0' \ar[r] \ar[d]^{t^{-1}} & N_1
\ar[r] \ar[d]^{t^{-1}} & N_2 \ar[r] \ar[d]^{t^{-1}} & 0 \\
0 \ar[r] & N_0' \ar[r] & N_1 \ar[r]  & N_2 \ar[r] & 0 }
$$

This gives an exact sequence of cokernels

$$
0 \longrightarrow G':=\frac{N'_0}{t^{-1}N'_0}\longrightarrow
G_I(M_1)\longrightarrow G_I(M_2)\longrightarrow 0 .$$ \noindent
Notice that
$$G'=\bigoplus_{n\in\N}\bigoplus_{i=0}^{\infty}\left [\frac{M_0\cap
I^iM_1}{M_0\cap I^{i+1}M_1}\right ]_n.$$ Let $u$ be another
indeterminate which we consider with degree one. Set
$$K_j:={\mathcal R}({\mathfrak
m}G_I(A),G_I(M_j))^+=\bigoplus_{s\in\Z}\bigoplus_{n\in\N}\bigoplus_{i=0}^
{\infty}\left [\frac{{\mathfrak
m}^sI^iM_j+I^{i+1}M_j}{I^{i+1}M_j}\right ]_n.$$ Notice that
$\frac{K_j}{\;u^{-1}K_j}=G_{\mathfrak m}(G_I(M_j)).$ Set
$$K'_0:=\ker(K_1\rightarrow K_2)=\bigoplus_{s\in\Z}\bigoplus_{n\in\N}
\bigoplus_{i=0}^{\infty}\left [\frac{M_0\cap ({\mathfrak
m}^sI^iM_1+I^{i+1}M_1)}{M_0\cap I^{i+1}M_1}\right ]_n.$$ We
consider the natural diagram
$$
\xymatrix { 0 \ar[r] & K_0'(1,0) \ar[r] \ar[d]^{u^{-1}} & K_1(1,0)
\ar[r] \ar[d]^{u^{-1}} & K_2(1,0) \ar[r] \ar[d]^{u^{-1}} & 0 \\
0 \ar[r] & K_0' \ar[r] & K_1 \ar[r]  & K_2 \ar[r] & 0 }
$$
\noindent which gives an exact sequence of cokernels

\begin{equation}\label{equation2}
0 \longrightarrow G'':=\frac{K'_0}{u^{-1}K'_0}\longrightarrow
G_{\mathfrak m}(G_I(M_1))\longrightarrow G_{\mathfrak
m}(G_I(M_2))\longrightarrow 0.
\end{equation}

Notice that
$$G''=\bigoplus_{s,n\in\N}\bigoplus_{i=0}^{\infty}\left
[\frac{M_0\cap ({\mathfrak m}^sI^iM_1+I^{i+1}M_1)}{M_0\cap
({\mathfrak m}^{s+1}I^iM_1+I^{i+1}M_1)}\right ]_n.$$ Let
$Q:=\mbox{Img}(K_0\longrightarrow K'_0)$ and let
$$P:=\ker(K_0\longrightarrow
K'_0)=\bigoplus_{s\in\Z}\bigoplus_{n\in\N}\bigoplus_{i=0}^{\infty}\left
[\frac{I^{i+1}M_0+{\mathfrak m}^sI^iM_0\cap
I^{i+1}M_1}{I^{i+1}M_0}\right ]_n.$$ We consider the natural
diagram
$$
\xymatrix { 0 \ar[r] & P(1,0) \ar[r] \ar[d]^{u^{-1}} & K_0(1,0)
\ar[r] \ar[d]^{u^{-1}} & Q(1,0) \ar[r] \ar[d]^{u^{-1}} & 0 \\
0 \ar[r] & P \ar[r] & K_0 \ar[r]  & Q \ar[r] & 0 }
$$
\noindent which gives an exact sequence of cokernels

\begin{equation}\label{equation3}
0 \longrightarrow P':=\frac{P}{u^{-1}P}\longrightarrow
G_{\mathfrak m}(G_I(M_0))\longrightarrow
Q':=\frac{Q}{u^{-1}Q}\longrightarrow 0.
\end{equation}

Notice that $$P'=\bigoplus_{s,n\in\N}\bigoplus_{i=0}^{\infty}\left
[\frac{I^{i+1}M_0+{\mathfrak m}^sI^iM_0\cap
I^{i+1}M_1}{I^{i+1}M_0+{\mathfrak m}^{s+1}I^iM_0\cap
I^{i+1}M_1}\right ]_n.$$ We consider the natural diagram
$$
\xymatrix { 0 \ar[r] & Q(1,0) \ar[r] \ar[d]^{u^{-1}} & K_0'(1,0)
\ar[r] \ar[d]^{u^{-1}} & L(1,0) \ar[r] \ar[d]^{u^{-1}} & 0 \\
0 \ar[r] & Q \ar[r] & K_0' \ar[r]  & L \ar[r] & 0 }
$$
\noindent where $L:=\mbox{coker}(Q\rightarrow K'_0).$ The
snake-lemma yields an exact sequence

\begin{equation}\label{equation4}
0 \longrightarrow W\longrightarrow Q'\longrightarrow
G''\longrightarrow V \longrightarrow 0
\end{equation}

\noindent where $W:=\ker(L(1,0)\rightarrow L)$ and
$V:=\mbox{coker}(L(1,0)\rightarrow L).$ We also have the exact
sequence

\begin{equation}\label{equation5}
0 \longrightarrow W\longrightarrow L(1,0)\longrightarrow
L\longrightarrow V \longrightarrow 0
\end{equation}

For $n\leq 1$ the coefficient modules of $u^n$ in $K_0={\mathcal
R}(\mathfrak m, G_I(M_0))^+$ and in $K_0'$ coincide, hence the
action of $u^{-1}$ on $L$ is nilpotent. Therefore the dimension of
$L$ is at most that of $G'',$ which is bounded by $D.$ Thus all
modules occurring in the exact sequence (\ref{equation5}) have
dimension at most $D.$

For any bigraded algebra $E=\bigoplus_{s,n\in\N}E_{s,n}$ consider
the Hilbert-Samuel functions $h_E(s,n):=\ell(E_{s,n})$ and its
Hilbert sums
$$h^{(1,0)}_E(s,n):=\sum_{v=0}^{s}h_E(v,n).$$

For $s,n\gg 0,$ we define the polynomial $h^D_{E}(s,n)$ of degree
$D-2$ as the Hilbert polynomial of $h_{E}(s,n)$ adding coefficient
zero to the terms of degree between $\dim E-2$ and $D-2.$

Using the additivity of the length function in (\ref{equation2}),
(\ref{equation3}) and (\ref{equation4}) leads to

{\footnotesize {\begin{equation}\label{equation6}
h_{(I,M_0)}^{*D\,(1,0)}(s,n)+h_{(I,M_2)}^{*D\,(1,0)}(s,n)-h_{(I,M_1)}^{*D\,
(1,0)}(s,n)
=h^{D\,(1,0)}_W(s,n)-h^{D\,(1,0)}_V(s,n)+h^{D\,(1,0)}_{P'}(s,n)
\end{equation}}}

Because of (\ref{equation5}) we have that
$$h^{D\,(1,0)}_W(s,n)-h^{D\,(1,0)}_V(s,n)=h^{D\,(1,0)}_L(s+1,n)-h^{D\,(1,0)}
_L(s,n)=h_L^{D}(s,n),$$ \noindent which, for $s,n\gg 0,$ is a
polynomial of degree at most $D-2$ because $L$ has dimension at
most $D.$

Because of (\ref{equation6}), for concluding the result we will
prove next that $h_{P'}^{D\,(1,0)}(s,n)$ is a polynomial of degree
at most $D-2,$ for $s,n\gg 0,$ or equivalently that
$h_{P'}^{(1,0)}(s,n)$ is a polynomial of degree at most $d-2,$ for
$s,n\gg 0.$ Notice that

$$h_{P'}^{(1,0)}(s,n)=\sum_{i=0}^{\infty}\ell \left[\frac{I^iM_0\cap I^{i+1}
M_1}{I^{i+1}M_0+{\mathfrak m}^{s+1}I^iM_0\cap I^{i+1}M_1}\right
]_n.$$ Set
$$S^i_{s,n}:=\left [\frac{I^iM_0\cap I^{i+1}M_1}{I^{i+1}M_0+{\mathfrak
m}^{s+1}I^iM_0\cap I^{i+1}M_1}\right ]_n,$$
$$C^i_{s,n}:=\left [\frac{I^iM_0}{{\mathfrak m}^{s+1}I^{i}M_0+I^iM_0\cap
I^{i+1}M_1}\right ]_n,$$ and  $$F^i_{s,n}:=\left [\frac{{\mathfrak
m}^{s}I^iM_0+I^iM_0\cap I^{i+1}M_1}{{\mathfrak
m}^{s+1}I^iM_0+I^{i+1}M_0}\right ]_n.$$

Consider the exact sequences

\begin{equation}\label{equation7}
0 \longrightarrow S^i_{s,n}\longrightarrow \left
[\frac{I^{i}M_0}{{\mathfrak m}^{s+1}I^{i}M_0+I^{i+1}M_0}\right
]_n\longrightarrow C^i_{s,n} \longrightarrow 0
\end{equation}

and

$$
0 \longrightarrow F^i_{s,n}\longrightarrow \left
[\frac{I^{i}M_0}{{\mathfrak m}^{s+1}I^{i}M_0+I^{i+1}M_0}\right
]_n\longrightarrow C^i_{s-1,n} \longrightarrow 0
$$
\noindent which yields

\begin{equation}\label{equation9}
h_{P'}^{(1,0)}(s,n)+ \sum_{i=0}^{\infty}\ell
(C^i_{s,n})-\sum_{i=0}^{\infty}\ell (C^i_{s-1,n})=
\sum_{i=0}^{\infty}\ell (F^i_{s,n}).
\end{equation}

Notice that, by equality (\ref{equation7}), for $s,n\gg 0,$
$\sum_{i=0}^{\infty}\ell (C^i_{s,n})$ is a polynomial of degree at
most $d-1.$ Hence, for $s,n\gg 0,$ $\sum_{i=0}^{\infty}\ell
(C^i_{s,n})-\sum_{i=0}^{\infty}\ell (C^i_{s-1,n})$ is a polynomial
of degree at most $d-2.$ On the other hand, if
$h_F(s,n):=\sum_{i=0}^{\infty}\ell (F^i_{s,n})$ then, it is clear
that $h^{(1,0)}_F(s,n)=h^{(1,0)}_{(I,M_0)}(s,n).$ Therefore, for
$s,n\gg 0,$ $h_F(s,n)$ is a polynomial of degree at most $d-2.$
Hence, by equality (\ref{equation9}), $h_{P'}^{(1,0)}(s,n)$ is a
polynomial of degree at most $d-2$, for $s,n\gg 0,$ as we claimed.
\end{proof}

In order to be able to formulate the main result in an efficient
way, we need a generalization of the notion of height of an ideal
to the case of modules (see \cite{Flenner-Manaresi}). We call the
number

$$\mbox{ht}_M(I):=\mbox{min }\{\dim M_{\mathfrak p}\mid {\mathfrak p}\in \mbox{ supp }M\cap V(I)\}$$
\noindent the $M$-height of $I$.

\begin{lem}\label{lemma1}
Let $(R, \mathfrak{m})$ be a local ring, $A$ a standard graded
Noetherian $R$-algebra and $M$ a finitely generated graded
$A$-module of dimension $D.$ Let $I\subseteq J$ be ideals of $A$
generated by linear forms such that $\mbox{ht}_M(I)>0$.  Then, for
all $k=0,\ldots, D-1$ we have that

\begin{itemize}
\item [(i)] $c^*_k(J,I^tM)=c^*_k(J,M)$ for all $t\in \N$;

\item [(ii)] $c^*_k(I,J^tM)=c^*_k(I,M)$ for all $t\in \N$;

\item [(iii)] $c^*_k(I,I^tM)=c^*_k(I,M)$ and
$c^*_k(J,J^tM)=c^*_k(J,M)$ for all $t\in \N$.
\end{itemize}
\end{lem}

\begin{proof}
In order to prove $(i)$ consider the exact sequence of graded
$A$-modules

$$
0 \longrightarrow I^tM\longrightarrow M\longrightarrow
\frac{M}{I^tM} \longrightarrow 0.
$$

From Theorem \ref{additivity} we have that
$$c^*_k(J,M)=c^*_{k,\,D}(J,I^tM)+c^*_{k,\,D}\left(J,\frac{M}{I^tM}\right).$$

Since $\mbox{ht}_M(I)>0$ we have that $\dim(M/I^tM)<D$ and $\dim
(I^tM)=D.$ Thus $c^*_{k,\,D}(J,I^tM)=c^*_k(J,I^tM)$ and
$c^*_{k,\,D}\left(J,\frac{M}{I^tM}\right)=0$, which proves $(i)$

The proof of $(ii)$ and $(iii)$ follows analogously.
\end{proof}

\begin{thm}\label{Rees-fraco}
Let $(R, \mathfrak{m})$ be a local ring, $A$ a standard graded
Noetherian $R$-algebra and $M$ a finitely generated graded
$A$-module of dimension $D.$ Let $I\subseteq J$ be ideals of $A$
generated by linear forms such that $\mbox{ht}_M(I)>0$.  If $I$ is
a reduction of $(J,M)$ then $c^*_k(I,M)=c^*_k(J,M)$ for all
$k=0,\ldots, D-1.$
\end{thm}

\begin{proof}
Since $I$ is a reduction of $(J,M)$ we have that
$I(J^rM)=J^{r+1}M=J(J^rM))$ for all $r\gg 0.$ Hence we have that
$G_{\mathfrak m}(G_I(J^rM))=G_{\mathfrak m}(G_J(J^rM))$ and thus
$$c^*_k(I, J^rM)=c^*_k(J, J^rM) \mbox{ for all } k=0,...,D-1.$$
Therefore the result follows by items $(ii)$ and $(iii)$ of Lemma
\ref{lemma1}.
\end{proof}

\section{$c^{\sharp}$-multiplicity sequence}

To prove a converse of Theorem \ref{Rees-fraco} we introduce
another multiplicity sequence, the $c^{\sharp}$-multiplicity
sequence, that is more suited for this purpose. The definition is
inspired by \cite{Ulrich-Validashti} and \cite{Achilles-Manaresi}.

In addition to the assumptions of the above section suppose that
$M$ is generated in degree zero. Again consider $G_I(M)$ as graded
by the internal grading.

\begin{defn}\label{c-sharp}
Let $D$ be any integer with $D\geq \dim M.$ We define the
\linebreak {\bf $c^{\sharp}$-multiplicity sequence} of $M$ with
respect to $I,$ as
$$c^{\sharp}_{k,\;D}(I,M):=c^D_k(A_1G_I(M)),\;k=0,\dots, D-1.$$
\noindent where $G_I(M)$ is graded by the internal grading. In the
case where $D=\dim M$ we simply write $c^{\sharp}_k(I,M)$ instead
of $c^{\sharp}_{k,\;\dim M}(I,M),\;k=0,\dots, \dim M-1.$
\end{defn}

To be more explicit, consider the standard bigraded $R$-algebra
\linebreak
$S^{\sharp}:=G_{\mathfrak{m}}(A_1G_I(A))=\oplus_{s,n=0}^{\infty}S^{\sharp}_
{s,n}$ with
$$S^{\sharp}_{s,n}=\oplus_{i=0}^{\infty}\left [\frac{{\mathfrak{m}}
^sI^iA_1+I^{i+1}}{\mathfrak{m}^{s+1}I^iA_1+I^{i+1}}\right ]_n=
\oplus_{i=0}^{n-1}\left
[\frac{\mathfrak{m}^sI^iA_1+I^{i+1}}{\mathfrak{m}^{s+1}I^iA_1+I^{i+1}}\right
]_n,$$ \noindent where $G_I(A)$ is graded by the internal grading,
and the finitely generated bigraded module over this algebra
$$T^{\sharp}=G_{\mathfrak{m}}(A_1G_I(M))=\oplus_{s,n=0}^{\infty}T^{\sharp}_
{s,n}$$ \noindent with
$$T^{\sharp}_{s,n}=\oplus_{i=0}^{\infty}\left [\frac{\mathfrak{m}^sI^iA_1M+I^
{i+1}M}{\mathfrak{m}^{s+1}I^iA_1M+I^{i+1}M}\right ]_n
=\oplus_{i=0}^{n-1}\left
[\frac{\mathfrak{m}^sI^iM+I^{i+1}M}{\mathfrak{m}^{s+1}I^iM+I^{i+1}M}\right
]_n,$$ \noindent where $G_I(M)$ is graded by the internal grading.

We denote the Hilbert-Samuel function $\ell_R(T^{\sharp}_{s,n})$
of $T^{\sharp}=G_{\mathfrak{m}}(A_1G_I(M))$ by
$h^{\sharp}_{(I,M)}(s,n)$ and its first Hilbert sum by $h^{\sharp
\,(1,0)}_{(I,M)}(s,n).$

Thus

$$h^{\sharp }_{(I,M)}(s,n)=\sum_{i=0}^{n-1}\ell_R \left [\frac{\mathfrak{m}
^sI^iM+I^{i+1}M}{\mathfrak{m}^{s+1}I^iM+I^{i+1}M}\right ]_n$$
\noindent and
$$h^{\sharp\,(1,0)}_{(I,M)}(s,n)=\sum_{i=0}^{n-1}\ell_R \left [\frac{I^iM}
{\mathfrak{m}^{s+1}I^iM+I^{i+1}M}\right ]_n.$$

For $s,n\gg 0,$ we define the polynomial $h^{\sharp\,
D}_{(I,M)}(s,n)$ of degree $D-2$ as the Hilbert polynomial of
$h^{\sharp}_{(I,M)}(s,n)$ adding coefficient zero to the terms of
degree between $\dim (A_1G_I(M))-2$ and $D-2.$

Thus, if $s,n\gg 0,$ the sequence $c^{\sharp}_{k,\; D}(I,M),
\;k=0,\ldots, D-1$ are the numerators of the leading coefficients
of the polynomial $h^{\sharp\, D\, (1,0)}_{(I,M)}(s,n).$

It will be useful to clarify the relationship between the two
multiplicity sequences $c^*$ and $c^{\sharp}.$

\begin{lem}\label{lemma2}
We use the same notation of Definition \ref{c-sharp}. Denote the
graded $G_I(A)$-module $G_I(M)/A_1G_I(M)$ by $B(I,M).$ Set
$b_k^D(I,M):=c_k^D(B(I,M)),$ $ k=0,\ldots, D-1.$ Then we have that

$$c^*_{k,\;D}(I,M)=c^{\sharp}_{k,\; D}(I,M)+b_k^D(I,M)$$
\end{lem}

\begin{proof}

Consider the exact sequence of $G_I(A)$-modules

$$0\longrightarrow A_1G_I(M)\longrightarrow G_I(M)\longrightarrow B(I,M)
\longrightarrow  0.$$ \noindent By the additivity of the
$c^D$-multiplicity sequence, Proposition \ref{proposition1}, we
have

$$c_k^D(G_I(M))=c_k^{D}(A_1G_I(M))+c_k^D(B(I,M)).$$
\noindent Recall that
$c_k^D(G_I(M))=c^*_{k,\;D}(I,M),\,c_k^{D}(A_1G_I(M))=c^{\sharp}_{k,\;
D}(I,M)$ and \linebreak $c_k^D(B(I,M))=b_k^D(I,M).$ Hence the
result follows.
\end{proof}

We will need the fact that the $c^{\sharp}$-multiplicity sequence
is additive on short exact sequences:

\begin{prop}({\it Additivity})\label{additivity2}
Let $(R, \mathfrak{m})$ be a local ring, $A$ a standard graded
Noetherian $R$-algebra, and $I$ an ideal of $A$ generated by
linear forms. If \linebreak $0 \longrightarrow M_0\longrightarrow
M_1\longrightarrow M_2\longrightarrow 0$ is an exact sequence of
finitely generated graded $A$-modules and $D$ an integer with
$D\geq\dim M_1.$ Then,
$$c^{\sharp}_{k,\;D}(I,M_1)=c^{\sharp}_{k,\;D}(I,M_0)+c^{\sharp}_{k,\;D}
(I,M_2)$$ \noindent for all $k=0,\ldots, D-1.$
\end{prop}

\begin{proof}
By Lemma \ref{lemma2} and Theorem \ref{additivity} it is enough to
show that the $b_k^D$-sequence is additive, that is
$$b_k^D(I,M_1)=b_k^D(I,M_0)+b_k^D(I,M_2)$$
\noindent for all $k=0,\ldots, D-1.$

Notice that $B(I,M_j)=\oplus_{n\in\N}[I^nM_j]_n,\,j=0,1,2.$ Set
$$G':=\ker (B(I,M_1)\longrightarrow B(I,M_2))=\oplus_{n\in\N}[M_0\cap
I^nM_1]_n.$$ \noindent We have the exact sequence

\begin{equation}\label{equation56}
0\longrightarrow G'\longrightarrow B(I,M_1)\longrightarrow
B(I,M_2)\longrightarrow  0.
\end{equation}

Set
$$L:=\mbox{coker} (B(I,M_0)\longrightarrow G')=\oplus_{n\in\N}\left [\frac
{M_0\cap I^nM_1}{I^nM_0}\right ]_n.$$ \noindent We have the exact
sequence

\begin{equation}\label{equation57}
0\longrightarrow B(I,M_0)\longrightarrow G'\longrightarrow
L\longrightarrow  0.
\end{equation}

Now (\ref{equation56}) and (\ref{equation57}) are exact sequences
of finitely generated graded modules of dimension at most $D.$
Hence we may compute the $c^D$-multiplicity sequence of graded
modules along these sequences. Using the additivity of this
multiplicity sequence as stated in Proposition \ref{proposition1}
we deduce that

$$b_k^D(I,M_1)=b_k^D(I,M_0)+b_k^D(I,M_2)+c_k^D(L).$$

To obtain that $c_k^D(L)=0$ we show that $L$ has dimension less
than $D.$ In fact, by Artin-Rees we have that

$$\frac{M_0\cap I^nM_1}{I^nM_0}=\frac{I^{n-c}(M_0\cap I^cM_1)}{I^nM_0}
\subseteq \frac{I^{n-c}M_0}{I^nM_0}.$$ \noindent Hence $\dim
(L)\leq \dim (P)$ where $P:=\oplus_{n\in\N}\left
[\frac{I^{n-c}M_0}{I^nM_0}\right ]_n.$ Clearly $P$ has dimension
less than the dimension of
$G_I(M_0)=\oplus_{n\in\N}\oplus_{i=0}^{\infty}\left
[\frac{I^{i}M_0}{I^{i+1}M_0}\right ]_n$ which is at most $D.$
\end{proof}

\begin{lem}\label{lemma3}
Let $(R, \mathfrak{m})$ be a local ring, $A$ a standard graded
Noetherian $R$-algebra and $M$ a finitely generated graded
$A$-module of dimension $D.$ Let $I\subseteq J$ be ideals of $A$
generated by linear forms such that $\mbox{ht}_M(I)>0$. Then, for
all $k=0,\ldots, D-1$ we have that

\begin{itemize}
\item [(i)] $c^{\sharp}_k(J,I^tM)=c^{\sharp}_k(J,M)$ for all $t\in
\N$;

\item [(ii)] $c^{\sharp}_k(I,J^tM)=c^{\sharp}_k(I,M)$ for all
$t\in \N$;

\item [(iii)] $c^{\sharp}_k(I,I^tM)=c^{\sharp}_k(I,M)$ and
$c^{\sharp}_k(J,J^tM)=c^{\sharp}_k(J,M)$ for all $t\in \N$.
\end{itemize}
\end{lem}

\begin{proof}
In order to prove $(i)$ consider the exact sequence of graded
$A$-modules

$$
0 \longrightarrow I^tM\longrightarrow M\longrightarrow
\frac{M}{I^tM} \longrightarrow 0.
$$

From Proposition \ref{additivity2} we have that
$$c^{\sharp}_k(J,M)=c^{\sharp}_{k,\,D}(J,I^tM)+c^{\sharp}_{k,\,D}\left
(J,\frac{M}{I^tM}\right).$$

Since $\mbox{ht}_M(I)>0$ we have that $\dim(M/I^tM)<D$ and $\dim
(I^tM)=D.$ Thus $c^{\sharp}_{k,\,D}(J,I^tM)=c^{\sharp}_k(J,I^tM)$
and $c^{\sharp}_{k,\,D}\left(J,\frac{M}{I^tM}\right)=0$, which
proves $(i)$

The proof of $(ii)$ and $(iii)$ follows analogously.
\end{proof}

\begin{rem}\label{rem}{\rm
Let $D=\dim M.$ If $\mbox{ht}_M(I)>0$ we have from the above Lemma
that, for all $k=0,\ldots, D-1$

\begin{itemize}
\item [(i)]
$c^{\sharp}_k(I,J^rM)=c^{\sharp}_k(I,I^rM)=c^{\sharp}_k(I,M)$ for
all $r\in \N$ and \item [(ii)]
$c^{\sharp}_k(J,I^sM)=c^{\sharp}_k(J,J^sM)=c^{\sharp}_k(J,M)$ for
all $s\in \N.$
\end{itemize}

That is, if
$$V^{(1,0,0,0)}(i,j,r,n):=\ell\left(\left [\frac{I^{j}J^{r}M}
{{\mathfrak{m}}^{i+1}I^{j}J^{r}M+I^{j+1}J^{r}M}\right
]_n\right),$$
 $$W^{(1,0,0,0)}(i,j,s,n):=\ell\left(\left [\frac{I^{s}J^{j}M}
{{\mathfrak{m}}^{i+1}I^{s}J^{j}M+I^{s}J^{j+1}M}\right
]_n\right),$$
$$H^{(1,0,0)}_{(I,M)}(i,j,n):=\ell\left(\left [\frac{I^{j}M}
{{\mathfrak{m}}^{i+1}I^{j}M+I^{j+1}M}\right ]_n\right),$$
\noindent and
$$H^{(1,0,0)}_{(J,M)}(i,j,n):=\ell\left(\left [\frac{J^{j}M}
{{\mathfrak{m}}^{i+1}J^{j}M+J^{j+1}M}\right ]_n\right)$$ \noindent
then for all $i,j,s,r,n\gg 0$ they becomes polynomials which
satisfy the following relations

\begin{equation}\label{equation101}
\overline{V^{(1,0,0,0)}(i,j,r,n)}=\overline{H^{(1,0,0)}_{(I,M)}(i,j+r,n)},
\end{equation}

\begin{equation}\label{equation112}
\overline{W^{(1,0,0,0)}(i,j,s,n)}=\overline{H^{(1,0,0)}_{(J,M)}(i,j+s,n)},
\end{equation}

\begin{equation}\label{equation102}
h^{\sharp\,(1,0)}_{(I,M)}(i,n)=H^{(1,1,0)}_{(I,M)}(i,n-1,n),
\end{equation}
\noindent and
\begin{equation}\label{equation103}
h^{\sharp\,(1,0)}_{(J,M)}(i,n)=H^{(1,1,0)}_{(J,M)}(i,n-1,n)
\end{equation}

\noindent where $\overline{F(i,j,r,n)}$ denotes, from now on, the
leading homogeneous part of the polynomial function $F(i,j,r,n).$
The equalities (\ref{equation101}) and (\ref{equation112}) follows
by $(i)$ and $(ii)$ respectively and the equalities
(\ref{equation102}) and (\ref{equation103}) follows by the
definition of $h^{\sharp\,(1,0)}_{(I,M)}(i,n)$ and
$h^{\sharp\,(1,0)}_{(J,M)}(i,n)$ respectively. By equations
(\ref{equation102}) and (\ref{equation103}) we have that, for all
$i,j,n\gg 0$, $H^{(1,0,0)}_{(J,M)}(i,j,n)$ and
$H^{(1,0,0)}_{(J,M)}(i,j,n)$ become polynomials of degree at most
$D-2.$ }
\end{rem}

\begin{thm}\label{Rees-fraco2}
Let $(R, \mathfrak{m})$ be a local ring, $A$ a standard graded
Noetherian $R$-algebra and $M$ a finitely generated graded
$A$-module of dimension $D.$ Let $I\subseteq J$ be ideals of $A$
generated by linear forms such that $\mbox{ht}_M(I)>0$.  If $I$ is
a reduction of $(J,M)$ then $c^{\sharp}_k(I,M)=c^{\sharp}_k(J,M)$
for all $k=0,\ldots, D-1.$
\end{thm}

\begin{proof}
Since $I$ is a reduction of $(J,M)$ we have that
$I(J^rM)=J^{r+1}M=J(J^rM))$ for all $r\gg 0.$ Hence we have that
$G_{\mathfrak m}(A_1G_I(J^rM))=G_{\mathfrak m}(A_1G_J(J^rM))$ and
thus
$$c^{\sharp}_k(I, J^rM)=c^{\sharp}_k(J, J^rM) \mbox{ for all } k=0,...,D-1.$$
Therefore the result follows by items $(ii)$ and $(iii)$ of Lemma
\ref{lemma3}.
\end{proof}

\section{Intertwining algebra and module}

In this section we recall the notions of interwining algebras and
intertwining modules introduced in \cite{Simis-Ulrich-Vasconcelos}
in the context of graded algebras which provide a strong criterium
for reductions of algebras. This algebras has been exploited on
several occasions by D. Kirby and D. Rees \cite{Kirby-Rees}, S.
Kleiman and A. Thorup \cite{Kleiman-Thorup} and D. Katz
\cite{Katz}. Their presentation can be immediately extended to the
version for modules we present here.

Let $(\mathcal R, \mathfrak m)$ be a Noetherian local ring, $A$ a
standard graded Noetherian $R$-algebra, $I\subseteq J$ ideals of
$A$ generated by linear forms and $M$ a finitely generated graded
$A$-module. Set
$$\mathcal A:={\mathcal R}(I,A)=\bigoplus_{i\in
\N}I^it^i\subseteq A\otimes_R R[t];$$ $$\mathcal B:={\mathcal
R}(J,A)=\bigoplus_{i\in \N}I^it^i\subseteq A\otimes_R R[t],$$
$${\mathcal R}(I,M):=\bigoplus_{i\in \N}I^iMt^i\subseteq M\otimes_R R[t],$$
\noindent and
$${\mathcal R}(J,M):=\bigoplus_{i\in \N}J^iMt^i\subseteq M\otimes_R
R[t].$$

Assigning degree zero to the variable $t$, the polynomial ring
\linebreak $A[t]:=A\otimes_R R[t]$ becomes a standard graded
Noetherian $R[t]$-algebra, and $M[t]:=M\otimes_R R[t]$ a finitely
generated module over this algebra. The Rees algebras $\mathcal A$
and $\mathcal B$ are homogeneous $R[t]$-subalgebras of $A[t],$ and
hence standard graded Noetherian  $R[t]$-algebras. Furthermore
${\mathcal R}(I,M)$ and ${\mathcal R}(J,M)$ are homogeneous
$\mathcal A$ and $\mathcal B$  submodules of $M[t]$ respectively,
thus they are finitely generated graded modules over $\mathcal A$
and $\mathcal B$ respectively. The grading so defined on this Rees
algebras and modules are also called {\it internal grading}.
Notice that with respect to this grading we have
$$[{\mathcal R}(I,M)]_n=\bigoplus_{i\in \N}[I^iM]_n,$$
\noindent and so on.

Let $u$ be a new variable which we also consider of degree zero
and set
$$C:={\mathcal R}(I,J):=\mathcal B\left [\bigoplus_{i\in \N}AI^iu^i\right ]
\subseteq A\otimes_R R[u,t]$$ \noindent Notice that
$A[u,t]:=A\otimes_R R[u,t]$ becomes a standard graded Noetherian
$R[u,t]$-algebra and $C$ a homogeneous $R[u,t]$-subalgebra of
$A[u,t]$, and hence a standard graded Noetherian $R[u,t]$-algebra.
Furthermore with respect to this grading, $C$ becomes a standard
graded Noetherian $R[u,t]$-algebra. This grading on $C$ is also
called {\it internal grading}. Notice that
$$C_n=\bigoplus_{i,j\in \N}[AI^iJ^j]_n.$$

Set
$$T_{J/I}(M):= \frac{C\left (\mathcal{R}(J,M)\right )}{C\left (\mathcal{R}
(I,M)\right )}.$$ \noindent where all the graded algebras and
graded modules involved are considered with the internal grading.
This grading on $T_{J/I}(M)$ is also called {\it internal
grading}. With this grading $T_{J/I}(M)$ becomes a finitely
generated Noetherian graded $C$-module which is called the {\it
intertwining module of} $I$ and $J$ with respect to $M.$ Notice
that
$$\left [T_{J/I}(M)\right ]_n=\bigoplus_{s,r\in \N}\left
[\frac{I^{s-1}J^{r+1}M} {I^sJ^rM}\right ]_n.$$ We say that $I$ is
a {\bf reduction} of $(J,M)$ if $IJ^nM=J^{n+1}M$ for at least one
positive integer $n.$

The following Theorem has been proved by A. Simis, B. Ulrich and
W.Vasconcelos \cite{Simis-Ulrich-Vasconcelos} in the context of
graded algebras (see also \cite[Theorem 1.153, p.
85]{Vasconcelos}). Their proof can be immediately extended to the
version for modules we present here.

\begin{thm}\label{interwining}
Let $(\mathcal R, \mathfrak m)$ be a Noetherian local ring, $A$ a
standard graded Noetherian $R$-algebra and $M$ a finitely
generated quasi-unmixed graded $A$-module. Let $I\subseteq J$ be
ideals of $A$ generated by linear forms such that
$\mbox{ht}_M(I)>0$. Then the following are equivalent:

\begin{itemize}
\item [(i)] $I$ is a reduction of $(J, M)$;

\item [(ii)] $\dim T_{J/I}(M)\leq \dim ({\mathcal R}(J,M))-1=\dim
M.$
\end{itemize}
\end{thm}

\begin{rem}{\rm
The implication $(i)\Rightarrow (ii)$ does not need the
quasi-unmixedness hypotheses for $M$ (see for example,
\cite[Proposition 1.149]{Vasconcelos}) this requirement is needed
only for the converse.}
\end{rem}

\begin{thm}\label{reesthm}
Let $(R, \mathfrak{m})$ be a Noetherian local ring, $A$ a standard
graded Noetherian $R$-algebra, $M$ a $D$-dimensional graded
$A$-module generated by finitely many homogeneous elements of
degree zero and $I\subseteq J$ ideals of $A$ generated by linear
forms such that $\mbox{ht}_M(I)>0$. Consider the following
statements:

\begin{enumerate}
\item [(i)] $I$ is a reduction of $(J,M)$;

\item [(ii)] $c^{\sharp}_k(I,M)=c^{\sharp}_k(J,M)$ for all
$k=0,...,D-1$.
\end{enumerate}
Then, $(i)$ implies $(ii)$ and if $M$ is quasi-unmixed the
converse also holds.
\end{thm}

\begin{proof}
The implication $(i)\Rightarrow (ii)$ has been proved in Theorem
\ref{Rees-fraco2}.

Conversely assume that $c^{\sharp}_k(I,M)=c^{\sharp}_k(J,M)$ for
all $k=0,...,D-1$. Notice that, by Theorem \ref{interwining}, it
is enough to prove that $\dim_{C(I,J)}\left(T_{J/I}(M)\right) \leq
\dim({\mathcal R}_{J}(M))-1=D.$ We will compute
$\dim_{C(I,J)}\left(T_{J/I}(M)\right)$ or, equivalently, ${\dim}
\left(G_{\mathfrak{m}}\left(T_{J/I}(M)\right)\right).$ Notice that

$$G:=G_{\mathfrak{m}}\left (T_{J/I}(M)\right )
=\bigoplus_{i,n\in \N}G_{i,n}$$ where
$$G_{i,n}:=\oplus_{j,r\in\N}\left [\frac{{\mathfrak{m}}^iI^{j-1}J^{r+1}M+I^jJ^rM}
{{\mathfrak{m}}^{i+1}I^{j-1}J^{r+1}M+I^jJ^rM}\right]_n.$$

Let $h_G(i,n)=\ell (G_{i,n})$ and let
$h^{(1,0)}_G(i,n)=\sum_{u=0}^{i}h_G(u,n)$ be its Hilbert-sum.
Notice that
$$h^{(1,0)}_G(i,n)=\sum_{j+r\leq n}\ell
\left(\left [\frac{I^{j-1}J^{r+1}M}
{{\mathfrak{m}}^{i+1}I^{j-1}J^{r+1}M+I^jJ^rM}\right ]_n\right).$$
For concluding the proof it is sufficient to show that, for
$i,n\gg 0,$ $h^{(1,0)}_G(i,n)$ is a polynomial of degree at most
$D-1.$ Set
$$F^{j,r}_{i,n}:=\left [\frac{I^{j-1}J^{r+1}M}
{{\mathfrak{m}}^{i+1}I^{j-1}J^{r+1}M+I^jJ^rM}\right ]_n,$$
$$P^{j,r}_{i,n}:=\left [\frac{{\mathfrak{m}}^{i+1}I^{j-1}J^{r+1}M\cap I^jJ^rM+I^jJ^{r+1}M}
{{\mathfrak{m}}^{i+1}I^{j}J^{r}M+I^{j+1}J^rM}\right ]_n$$ and
$$Q^{j,r}_{i,n}:=\left [\frac{I^{j}J^{r}M}
{{\mathfrak{m}}^{i+1}I^{j}J^{r}M+I^{j+1}J^rM}\right ]_n.$$

Notice that

$$h^{(1,0)}_G(i,n)=\sum_{j+r\leq n}\ell \left (\left [F^{j,r}_{i,n} \right ]_n\right ).$$
\noindent Therefore to conclude the proof it is sufficient to show
that, for $i,j,r,n\gg 0,$ $\ell \left (\left [F^{j,r}_{i,n} \right
]_n\right )$ becomes a polynomial of degree at most $D-3.$

Consider the exact sequence

$$
0 \longrightarrow P^{j,r}_{i,n}\longrightarrow
Q^{j,r}_{i,n}\longrightarrow Q^{j-1,r+1}_{i,n}\longrightarrow
F^{j,r}_{i,n} \longrightarrow 0
$$
\noindent which yields

\begin{equation}\label{equation114}
 \ell \left (\left [F^{j,r}_{i,n} \right
]_n\right )=\ell (Q^{j-1,r+1}_{i,n})-\ell (Q^{j,r}_{i,n})+\ell
(P^{j,r}_{i,n})
\end{equation}

Notice that, by equality (\ref{equation101}) of Remark \ref{rem},
 $\overline{\ell (Q^{j,r}_{i,n})}=\overline{H^{(1,0,0)}_{(I,M)}(i,j+r,n)}$ for all
$i,j,r,n\gg 0.$ Hence, for all $i,j,r,n\gg 0,$ $\ell
(Q^{j-1,r+1}_{i,n})-\ell (Q^{j,r}_{i,n})$ is a polynomial of
degree at most $D-3.$ Therefore, by equality (\ref{equation114}),
it remains to prove that,  for $i,j,r,n\gg 0,$ $\ell
(P^{j,r}_{i,n})$ is polynomial of degree at most $D-3.$

We observe that, by Artin-Rees,
$$
P^{j,r}_{i,n}\subseteq N^{j,r}_{i,n}:=\left
[\frac{{\mathfrak{m}}^{i+1-c}I^{j}J^{r}M+I^jJ^{r+1}M}
{{\mathfrak{m}}^{i+1}I^{j}J^{r}M+I^{j+1}J^rM}\right ]_n.
$$
Thus, it is enough to prove that, for $i,j,r,n\gg 0,$ $\ell
(N^{j,r}_{i,n})$ is a polynomial of degree at most $D-3.$

Consider now the exact sequences

{\footnotesize{\begin{equation}\label{equation116} 0
\longrightarrow \left
[\frac{{\mathfrak{m}}^{i+1-c}I^{j}J^{r}M+I^{j+1}J^{r}M}
{{\mathfrak{m}}^{i+1}I^{j}J^{r}M+I^{j+1}J^rM}\right
]_n\longrightarrow N^{j,r}_{i,n}\longrightarrow \left
[\frac{{\mathfrak{m}}^{i+1-c}I^{j}J^{r}M+I^{j}J^{r+1}M}
{{\mathfrak{m}}^{i+1-c}I^{j}J^{r}M+I^{j+1}J^rM}\right ]_n
\longrightarrow 0,
\end{equation}}}

{\footnotesize{\begin{equation}\label{equation117} 0
\longrightarrow \left
[\frac{{\mathfrak{m}}^{i+1-c}I^{j}J^{r}M\!+\!I^{j}J^{r+1}M}
{{\mathfrak{m}}^{i+1-c}I^{j}J^{r}M\!+\!I^{j+1}J^rM}\right
]_n\longrightarrow Q^{j,r}_{i-c,n}\longrightarrow \left
[\frac{I^{j}J^{r}M}
{{\mathfrak{m}}^{i+1-c}I^{j}J^{r}M\!+\!I^{j}J^{r+1}M}\right ]_n
\longrightarrow 0
\end{equation}}}

\noindent and

\begin{equation}\label{equation118} 0 \longrightarrow
\left [\frac{{\mathfrak{m}}^{i+1-c}I^{j}J^{r}M+I^{j+1}J^{r}M}
{{\mathfrak{m}}^{i+1}I^{j}J^{r}M+I^{j+1}J^rM}\right
]_n\longrightarrow Q^{j,r}_{i,n}\longrightarrow
Q^{j,r}_{i-c,n}\longrightarrow 0
\end{equation}

\noindent which yields

$$\begin{array}{lll}
\vspace{0.3cm}
 \ell (N^{j,r}_{i,n}) & =    & \ell\left(\left [\frac{{\mathfrak{m}}^{i+1-c}I^{j}J^{r}M+I^{j}J^{r+1}M}
{{\mathfrak{m}}^{i+1-c}I^{j}J^{r}M+I^{j+1}J^rM}\right ]_n\right)+
\ell \left(\left
[\frac{{\mathfrak{m}}^{i+1-c}I^{j}J^{r}M+I^{j+1}J^{r}M}
{{\mathfrak{m}}^{i+1}I^{j}J^{r}M+I^{j+1}J^rM}\right ]_n\right) \\
\vspace{0.3cm} & = & \ell (Q^{j,r}_{i,n})-\ell \left(\left
[\frac{I^{j}J^{r}M}
{{\mathfrak{m}}^{i+1-c}I^{j}J^{r}M\!+\!I^{j}J^{r+1}M}\right
]_n\right)
\end{array}$$
\noindent where the first equality follows by (\ref{equation116})
and the last equality follows by (\ref{equation117}) and
(\ref{equation118}).

Therefore, for $i,j,r,n\gg 0,$

\begin{equation}\label{equation119}
\overline{\ell (N^{j,r}_{i,n})}=\overline{\ell (Q^{j,r}_{i,n})}-
\overline{\ell \left(\left [\frac{I^{j}J^{r}M}
{{\mathfrak{m}}^{i+1-c}I^{j}J^{r}M\!+\!I^{j}J^{r+1}M}\right
]_n\right)}
\end{equation}

Notice that, by equality (\ref{equation101}) of Remark \ref{rem}
 $\overline{\ell (Q^{j,r}_{i,n})}=\overline{H^{(1,0,0)}_{(I,M)}(i,j+r,n)}$ for all
$i,j,r,n\gg 0$ and by equality (\ref{equation112}) of Remark
\ref{rem}

$$\overline{\ell \left(\left [\frac{I^{j}J^{r}M}
{{\mathfrak{m}}^{i+1-c}I^{j}J^{r}M\!+\!I^{j}J^{r+1}M}\right
]_n\right)} =\overline{H^{(1,0,0)}_{(J,M)}(i,j+r,n)}$$

Hence, for $i,j,r,n\gg 0,$
$$\overline{\ell (N^{j,r}_{i,n})}=\overline{H^{(1,0,0)}_{(I,M)}(i,j+r,n)}
-\overline{H^{(1,0,0)}_{(J,M)}(i,j+r,n)}$$

Notice that the leading coefficients of the Hilbert-Samuel
polynomials of $H^{(1,0,0)}_{(I,M)}(i,j+r,n)$ and
$H^{(1,0,0)}_{(J,M)}(i,j+r,n)$ are described as specific sums of
$c^{\sharp}_k(I,M)$ and $c^{\sharp}_k(J,M),$ with $k=0,...,D-1,$
respectively, which by assumption must coincide. Hence, for
$i,j,r,n\gg 0,$ $\ell (N^{j,r}_{i,n})$ is a polynomial of degree
at most $D-3$ as we claimed.

\end{proof}

\section{Multiplicity sequence for arbitrary modules}

We are now ready to introduce the main object of this paper, the
multiplicity sequence of a module. Here the ideal of the previous
section will be replaced by a module $E.$

Let $(R, \mathfrak{m})$ be a Noetherian local ring, $E$ a
submodule of the free $R$-module $R^p,$ and $N$ a finitely
generated $R$-module of dimension $d.$  The symmetric algebra
$A:=\mbox{Sym}(R^p)=\oplus S_n(R^p)$ of $R^p$ is a polynomial ring
$R[T_1,\ldots ,T_p].$ If $h=(h_1,\ldots ,h_p)\in R^p,$ then we
define the element $w(h)=h_1T_1 + \ldots +h_pT_p\in A.$ We denote
by ${\cal R}(E):=\oplus {\cal R}_n(E)$ the subalgebra of $A$
generated in degree one by $\{w(h): h\in E\}$ and call it {\it the
Rees algebra} of $E$. Then ${\cal R}(E)$ has dimension $d+p$.
Consider the $A$-ideal $I$ generated by ${\cal R}_1(E)$ and the
$A$-module $M:=A\otimes_R N.$ Notice That $I$ is an $A$-ideal
generated by linear forms and $M$ is a finitely generated graded
$A$-module of dimension $d+p$ that is generated in degree zero.

\begin{defn}\label{BR-multiplicity}
We define the {\bf multiplicity sequence} associated to the module
$E$ with respect to $N$ by
$$c_k(E,N):=c^{\sharp}_k(I,M),\;\;\;\;k=0,\dots, d+p-1.$$
\end{defn}

To be more explicit,
$$\left [\frac{I^iM}{\mathfrak{m}^{s+1}I^iM+I^{i+1}M}\right ]_n=
\frac{{\cal R}_i(E)S_{n-i}(R^p)N}{\mathfrak{m}^{s+1}{\cal
R}_i(E)S_{n-i}(R^p)N+{\cal R}_{i+1}(E)S_{n-i-1}(R^p)N}$$ \noindent
for $0\leq i\leq n-1.$ Thus the Hilbert function of
$T^{\sharp}=G_{\mathfrak{m}}(A_1G_I(M))$ is

$$h^{\sharp\,(1,0)}_{(I,M)}(s,n)=\sum_{i=0}^{n-1}\ell_R \left[\frac{{\cal R}
_i(E)S_{n-i}(R^p)N}{\mathfrak{m}^{s+1}{\cal
R}_i(E)S_{n-i}(R^p)N+{\cal R}_{i+1}(E)S_{n-i-1}(R^p)N}\right ],$$
\noindent which, for $s,n\gg 0,$ becomes a polynomial of degree at
most $d+p-1$ whose leading coefficients are
$c_k(E,N),\;k=0,\ldots,d+p-1.$ If $N=R$ we simply write $c_k(E)$
instead of $c_k(E,N)$ for $k=0,\ldots, d+p-1.$

\begin{rem}{\rm
If $E$ has finite colength in $R^p$ then, for $s\gg 0,$ we have
that $\mathfrak{m}^{s+1}{\cal R}_i(E)S_{n-i}(R^p)\subseteq {\cal
R}_{i+1}(E)S_{n-i-1}(R^p).$ Hence, in this context

$$
\begin{array}{lll}
\vspace{0.3cm} h^{\sharp\,(1,0)}_{(I,M)}(s,n) & = &
\sum_{i=0}^{n-1}\ell_R \left[\frac{{\cal
R}_i(E)S_{n-i}(R^p)}{\mathfrak{m}^{s+1}{\cal
R}_i(E)S_{n-i}(R^p)+{\cal R}_{i+1}(E)S_{n-i-1}(R^p)}\right ]
\\
\vspace{0.3cm} & = & \sum_{i=0}^{n-1}\ell_R \left[\frac{{\cal
R}_i(E)S_{n-i}(R^p)}{{\cal R}_{i+1}(E)S_{n-i-1}(R^p)}\right ] \\
\vspace{0.3cm} & = & \ell_R \left[\frac{S_{n}(R^p)}{{\cal
R}_{n}(E)}\right ]
\end{array}
$$
\noindent Thus in this case $c_0(E)=e_{BR}(E)$ and $c_k(E)=0$ for
all $k=1,\ldots, d+p-1.$ In case that $E$ is an ideal $J$ of $R$
then  the Buchsbaum-Rim multiplicity sequence $c_k(E,N)$ coincides
with the Achilles-Manaresi multiplicity sequence $c_k(J, N)$ for
all $k=0,\ldots, d.$ }
\end{rem}

Theorem \ref{reesthm} immediately gives the following result:

\begin{thm}\label{reesmodule} Let $(R, \mathfrak{m})$ be a
Noetherian local ring, $E\subseteq F\subseteq R^p$ be $R$-modules
and write $I:={\cal R}_1(E)$ for the corresponding ideal of
$A:=\mbox{Sym}(R^p).$ Let $N$ be a $d$-dimensional finitely
generated $R$-module and set $M:=A\otimes_R N.$ Assume that
$\mbox{ht }_M(I)>0.$ Consider the following statements:

\begin{enumerate}
\item [(i)] $E$ is a reduction of $(F,N)$;

\item [(ii)] $c_k(E,N)=c_k(F,N)$ for all $k=0,...,d+p-1$.
\end{enumerate}
Then, $(i)$ implies $(ii)$ and if $N$ is quasi-unmixed the
converse also holds.
\end{thm}

\end{document}